\documentclass[12pt]{article}
\usepackage{amssymb}
\usepackage{amsmath}
\usepackage{amsbsy}
\newtheorem{lemma}{Lemma}[section]

\newtheorem{theorem}{Theorem}[section]

\def\proclaim#1{\par \bigskip\noindent {\bf #1}\bgroup\it\ }
\def\endproclaim{\egroup\par\bigskip}

\newbox\TempBox \newbox\TempBoxA
\newcommand{\non}{\nonumber \\}
\def\pr{\textsf{P}} 
\def\ep{\textsf{E}} 
\def\Cov{\textsf{Cov}} 

\def\text#1{\mbox{\rm #1}}

\def\underwiggle 1{
\ifmmode\setbox\TempBox=\hbox{$ 1$}\else\setbox\TempBox=\hbox{
1}\fi \setbox\TempBoxA=\hbox to \wd\TempBox{\hss\char'176\hss}
\rlap{\copy\TempBox}\smash{\lower9pt\hbox{\copy\TempBoxA}} }

\begin{document}

\begin{center}
{\Large On maxima of periodograms  of stationary processes\footnote{ Research supported by National Natural Science Foundations of China
(10571159, 10671176) and Specialized Research Fund for the Doctor Program of Higher Education (20060335032). }}
\end{center}

\begin{center}
{{Zhengyan Lin\footnote{Email: zlin@zju.edu.cn} and Weidong Liu\footnote{Email: liuweidong99@gmail.com. }  }}

{{\it Department of Mathematics, Zhejiang University, Hangzhou 310027, China } }

\end{center}

\noindent {\rm {\bf Abstract.} We consider the limit distribution of maxima of periodograms for stationary processes. Our method is based on
$m$-dependent approximation for stationary processes and a moderate deviation result.

 \noindent {\bf Keywords:} stationary process,  periodogram, $m$-dependent
approximation.

\bigskip
\noindent{\bf AMS 2000 subject classification:} 62M15; 60F05 }
\section{Introduction}
 \setcounter{equation}{0}
Let $\{\varepsilon_{n}; n\in Z\}$ be independent and identically distributed (i.i.d.) random variables and $g$ be a measurable function such that
 \begin{equation}\label{a0}
 X_{n}=g(\cdots,\varepsilon_{n-1},\varepsilon_{n})
 \end{equation}
is a well-defined random variable. Then $\{X_n; n\in Z\}$ presents a huge class of processes. In particular, it contains the linear process and
nonlinear processes including the threshold AR (TAR) models, ARCH models, random coefficient AR (RCA) models, exponential AR (EAR) models and so
on. Wu and Shao \cite{wu3} argued that many nonlinear time series are stationary causal with one-sided representation (\ref{a0}). Let
\begin{equation*}
I_{n,X}(\omega)=n^{-1}\Big{|}\sum_{k=1}^{n}X_{k}\exp(\text{i}\omega k)\Big{|}^{2},\quad \omega\in [0, \pi],
\end{equation*}
be the periodogram of random variables $X_{1},\cdots,X_{n}$ and denote
\begin{equation*}
M_{n}(X)=\max_{1\leq j\leq q}I_{n,X}(\omega_{j}),\quad \omega_{j}=2\pi j/n,
\end{equation*}
where $q=q_{n}=\max\{j:0<\omega_{j}<\pi\}$ so that $q\sim n/2$.

If $X_{1}, X_{2},\cdots$ are i.i.d. random variables with $N(0, 1)$ distribution, then $\{I_{n,X}(\omega_{j}); 1\leq j\leq q\}$ is a sequence of
i.i.d. standard exponential random variables. It is well-known that (cf. Brockwell and Davis \cite{brock})
\begin{equation}\label{a1}
M_{n}(X)-\log q\Rightarrow G,
\end{equation}
 where $\Rightarrow$ means convergence in distribution and $G$ has the
standard Gumbel distribution $\Lambda(x)=\exp(-\exp(-x))$, $x\in R$. However, in the non-Gaussian case, the independence of $I_{n,X}(\omega_{j})$
is not guaranteed in general, and therefore (\ref{a1}) is not trivial. When $X_{1}, X_{2},\cdots$ are i.i.d. random variables, Davis and Mikosch
\cite{davis}  established (\ref{a1}) with the assumptions that $\ep X_{1}=0$, $\ep X^{2}_{1}=1$ and $\ep |X_{1}|^{s}<\infty$ for some $s>2$. They
also conjectured that the condition $\ep X^{2}_{1}\log^{+}|X_{1}|<\infty$ is
 sufficient for (\ref{a1}). Moreover, a similar result was established in their paper for the two-sided linear process
$X_{n}=\sum_{j\in Z}a_{j}\varepsilon_{n-j}$ under the conditions that $\ep |\varepsilon_{0}|^{s}<\infty$ for some $s>2$ and
\begin{equation}\label{a2}
\sum_{j\in Z}|j|^{1/2}|a_{j}|<\infty.
\end{equation}
The key step in Davis and Mikosch \cite{davis} is the following approximation (cf. Walker \cite{walker})
\begin{equation}\label{dm1}
\max_{\omega\in[0,\pi]}\Big{|}\frac{I_{n,X}(\omega)}{2\pi f(\omega)}-I_{n,\varepsilon}(\omega)\Big{|}\rightarrow_{\pr} 0.
\end{equation}

Generally, it is very difficult to check (\ref{dm1}) for the stationary process defined in (\ref{a0}). In this paper, we shall establish
(\ref{a1}) (or an analogous result) for (\ref{a0}) under some regularity conditions. Let us take a look at the linear process first. In this case,
$X_{n}=\sum_{j=-m}^{m}a_{j}\varepsilon_{n-j}+\sum_{|j|>m}a_{j}\varepsilon_{n-j}$, $m>0$. Under the assumptions of $\sum_{j\in Z}|a_{j}|<\infty$
and $\ep|\varepsilon_{0}|<\infty$, $\sum_{|j|>m}a_{j}\varepsilon_{n-j}\rightarrow 0$ in probability as $m\rightarrow\infty$. This implies that the
linear process behaves like a process which is block-wise independent. In fact, many time series, such as the GARCH model, have such property.
Such an analysis suggests us to approximate $X_{n}$ by $\ep[X_{n}|\varepsilon_{n-m},\cdots,\varepsilon_{n}]$. This method has been employed in
Hsing and Wu \cite{hsing} to establish the asymptotic normality of a weighted $U$-statistic.

By the $m$-dependent approximation developed in Section 3, we show that, for proving (\ref{a1}), the condition (\ref{a2}) can be weakened to
$\sum_{|j|\geq n}|a_{j}|=o(1/\log n)$. Meanwhile, the moment condition on $\varepsilon_{0}$ can also be weakened to
$\ep\varepsilon^{2}_{0}I\{|\varepsilon_{0}|\geq n\}=o(1/\log n)$. This in turn proves that
 the conjecture by Davis and Mikosch \cite{davis} is true.
Furthermore, it is shown that (\ref{a1}) still holds for the general process defined in (\ref{a0}).

Below, we explain how (\ref{a1}) (or the analogous result) can be used for detecting periodic components in a time series (see also Priestley
\cite{priestley}). Let us consider the model
\begin{equation*}
Z_{t}=\mu+S(t)+X_{t}\quad t=1,2,\ldots, n,
\end{equation*}
where $X_{t}$ is a stationary time series with mean zero and the deterministic part
\begin{equation*}
S(t)=A_{1}\cos(\gamma_{1}t+\phi_{1})
\end{equation*}
is a sinusoidal wave at frequency $\gamma_{1}\neq 0$ with the amplitude $A_{1}\neq 0$ and the phase $\phi_{1}$. Without loss of generality, we
assume $\mu=0$.
 A test statistic for the null hypothesis $H_{0}: S(t)\equiv 0$ against the
 alternative  $H_{1}: S(t)=A_{1}\cos(\gamma_{1}t+\phi_{1})$
  is
\begin{equation}\label{d1}
g_{n}(Z)=\frac{\max_{1\leq i\leq q}I_{n,Z}(\omega_{i})/\hat{f}(\omega_{i})}{\sum_{i=1}^{q}I_{n,Z}(\omega_{i})/\hat{f}(\omega_{i})},
\end{equation}
where $\hat{f}(\omega)$ is an estimator of $f(\omega)$, the spectral density of $Z_{t}$.
 This statistic was proposed by
Fisher \cite{fisher}, who assumed  that  $X_{t}$ is a white Gaussian series and thus chose $\hat{f}(\omega)\equiv 1$. Often, however, it is not
reasonable, as a null hypothesis, to assert that the observations are independent. Hence, Hannan \cite{hannan1} assumed that $X_{t}=\sum_{j\in
Z}a_{j}\varepsilon_{t-j}$ with $\varepsilon_{t}$ being i.i.d. normal and $\{a_{j}\}$ satisfying some conditions. The results in Section 2 make it
possible to obtain the asymptotic distribution of $g_{n}(Z)$ under $H_{0}$, for a class of general processes rather than the linear process, and
without the requirement of the normality for $\varepsilon_{t}$; see Remark 2.4 for more details.

Sometimes we might suspect that the series might contain several periodic components. In this case, we should test $H_{0}: S(t)\equiv 0$ against
the
 alternative $H_{1}:
 S(t)=\sum_{k=1}^{r}A_{k}\cos(\gamma_{k}t+\phi_{k})$, where $r (>1)$ is
 the possible number of peaks. Assuming that $X_{t}$ is a white Gaussian series, Shimshoni \cite{shim} and Lewis and Fieller \cite{lewis}
 proposed the statistic
 \begin{equation*}
 U_{Z}(r)=\frac{I_{n,q-r+1}(Z)}{\sum_{i=1}^{q}I_{n,Z}(\omega_{i})}
 \end{equation*}
 for
 detecting $r$ peaks. Here $I_{n,1}(Z)\leq I_{n,2}(Z)\leq\cdots\leq I_{n,q}(Z)$ are
 the order statistics of the periodogram ordinates
 $I_{n,Z}(\omega_{i})$, $1\leq i\leq q$. The exact (and asymptotic) null distribution of $U_{Z}(r)$
  can be found in Hannan \cite{hannan2} and Chiu \cite{chiu}. In the
 latter paper, the test statistic
 $R_{Z}(\beta)=I_{n,q}(Z)/\sum_{j=1}^{[q\beta]}I_{n,j}(Z)$, $0<\beta<1$, was
 given. Our results may be useful for obtaining the asymptotic
 distribution of $R_{Z}(\beta)$ when $X_{n}$ is defined in
 (\ref{a0}).

 The paper is
organized as follows. Our main results Theorems 2.1 and 2.2 will be presented in Section 2. In Section 3, we develop the $m$-dependent
approximation for the Fourier transforms of stationary processes. The proofs of main results will be given in Sections 4 and 5. Throughout the
paper, we let $C$, $C_{(\cdot)}$ denote
 positive constants and their values may be different in different contexts. When
$\delta$ appears, it usually means every $\delta>0$ and may be different in every place. For two real sequences $\{a_{n}\}$ and $\{b_{n}\}$, write
$a_{n}=O(b_{n})$ if there exists a constant $C$ such that $|a_{n}|\leq C|b_{n}|$ holds for large $n$, $a_{n}=o(b_{n})$ if
$\lim_{n\rightarrow\infty}a_{n}/b_{n}=0$ and $a_{n}\asymp b_{n}$ if $C_{1}b_{n}\leq a_{n}\leq C_{2}b_{n}$. With no confusion, we let $|\cdot|$
denote the $d$-dimensional Euclidean norm ($d\geq 1$) or the norm of a $d\times d$ matrix A, defined by $|\text{A}|=\max_{|x|\leq 1, x\in
R^{d}}|\text{A}x|$.
 \section{Main results}
 \setcounter{equation}{0}

We first consider the two-sided linear process. Let
\begin{equation}\label{m1}
Y_{n}=\sum_{j\in Z}a_{j}\varepsilon_{n-j}, \mbox{~and~} X_{n}=h(Y_{n})-\ep h(Y_{n}),
\end{equation}
 where $\sum_{j\in Z}|a_{j}|<\infty$ and
$h$ is a Lipschitz continuous function. Let us redefine
$$I_{n,1}(X)\leq I_{n,2}(X)\leq\cdots\leq I_{n,q}(X)$$ as
 the order statistics of the periodogram ordinates
 $I_{n,X}(\omega_{j})/(2\pi f(\omega_{j}))$, $1\leq j\leq q$, where $f(\omega)$ is the spectral density function of
 $\{X_{n}\}$, defined by
\begin{equation*}
f(\omega)=\frac{1}{2\pi}\sum_{k\in Z}\ep X_{0}X_{k}\exp(\text{i}k\omega)
\end{equation*}
 and
 satisfies
\begin{equation}\label{c1}
f^{*}:=\min_{\omega\in R}f(\omega)>0.
\end{equation}
  Note that $f(\omega)\equiv \ep X^{2}_{1}/(2\pi)$ if
 $X_{1}, X_{2},\cdots$ are i.i.d. centered random variables.

\begin{theorem} Let $X_{n}$ be defined in (\ref{m1}). Suppose that (\ref{c1}) holds, and
\begin{equation}\label{a3}
 \ep \varepsilon_{0}=0, ~\ep
\varepsilon^{2}_{0}=1\mbox{ and }\sum_{|j|\geq n}|a_{j}|=o(1/\log n).
\end{equation}
 (i). Suppose that $h(x)=x$, and
\begin{equation}\label{a4}
\ep\varepsilon^{2}_{0}I\{|\varepsilon_{0}|\geq n\}=o(1/\log n).
\end{equation}
Then
\begin{equation}\label{a5}
I_{n,q}(X)-\log q\Rightarrow G,
\end{equation}
where $G$ has the
standard Gumbel distribution $\Lambda(x)=\exp(-\exp(-x))$, $x\in R$.\\
 (ii). Suppose  $h$
is a Lipschitz continuous function on $R$. If (\ref{a4}) is strengthened to $\ep\varepsilon^{2}_{0}I\{|\varepsilon_{0}|\geq n\}=o(1/(\log
n)^{2})$, then (\ref{a5}) holds.
\end{theorem}

{\bf Remark 2.1.} From Theorem 2.1, we derive the asymptotic distribution of the maximum of the periodogram. Note that (\ref{a4}) is implied by
$\ep \varepsilon^{2}_{0}\log^{+} |\varepsilon_{0}|<\infty$. Hence the conjecture in Davis and Mikosch \cite{davis} is true. In order to show
$\max_{1\leq j\leq q}I_{n,X}(\omega_{j})/(2\pi f(\omega_{j}))-\log q\Rightarrow G$ when $X_{n}=\sum_{j\in Z}a_{j}\varepsilon_{n-j}$,  Davis and
Mikosch \cite{davis} used the following approximation
\begin{equation}\label{dm}
\max_{\omega\in[0,\pi]}\Big{|}\frac{I_{n,X}(\omega)}{2\pi f(\omega)}-I_{n,\varepsilon}(\omega)\Big{|}\rightarrow_{\pr} 0
\end{equation}
which requires the condition (\ref{a2}). Obviously, our condition in (\ref{a3}) is weaker than (\ref{a2}). They also required $\ep
|\varepsilon_{0}|^{s}<\infty$ for some $s>2$, which is stronger than (\ref{a4}). Moreover, it is difficult to prove (\ref{dm}) for the nonlinear
transforms of linear processes considered in {\em (ii)}.

{\bf Remark 2.2.} The (weak) law of logarithm for the maximum of the periodogram is a simple consequence of Theorem 2.1. Under conditions on the
smoothness of the characteristic function of $\varepsilon_{n}$, An et al. \cite{an} proved the (a.s.) law of logarithm for the maximum of the
periodogram.

In the following, we will give a theorem when $X_{n}$ satisfies the general form in (\ref{a0}). Of course, we should impose some dependency
conditions on $X_{n}$. For the reader's convenience, we list the following notations.
\begin{itemize}
\item{$\mathcal{F}_{i,j}:=(\varepsilon_{i},\cdots,\varepsilon_{j})$, $-\infty\leq i\leq j\leq \infty$.} \item{$Z\in L^{p}$ if
$\|Z\|_{p}:=(\ep|Z|^{p})^{1/p}<\infty$.   } \item{$\{\varepsilon^{*}_{i}, i\in Z\}$ is an independent copy of $\{\varepsilon_{i}, i\in Z\}$.}
\item{$\theta_{n,p}:=\|X_{n}-X^{*}_{n}\|_{p}$, where $X^{*}_{n}=g(\cdots,\varepsilon_{-1},\varepsilon^{*}_{0},\mathcal{F}_{1,n})$.}
\item{$\Theta_{n,p}:=\sum_{i\geq n}\theta_{i,p}$.}
\end{itemize}

{\bf Remark 2.3.}  $\theta_{n,p}$ is called the physical dependence measure by Wu \cite{wu}. An advantage of such dependence measure is that it is
easily verifiable.

\begin{theorem}Let $X_{n}$ be defined in (\ref{a0}) and (\ref{c1}) holds. Suppose that
$\ep X_{0}=0$, $\ep |X_{0}|^{s}<\infty$ for some $s>2$ and $\Theta_{n,s}=o(1/\log n)$.
 Then (\ref{a5}) holds.
\end{theorem}

{\bf Remark 2.4.} To derive the asymptotic distribution (under $H_{0}$) of $g_{n}(Z)$ defined in (\ref{d1}) from Theorem 2.2, we should prove
\begin{equation}\label{re1}
|q^{-1}\sum_{i=1}^{q}I_{n,Z}(\omega_{i})/(2\pi f(\omega_{i}))-1|=o_{\pr}(1/\log n)
\end{equation}
and choose $\hat{f}(\omega)$, an estimator of $f(\omega)$, to satisfy
\begin{equation}\label{re2}
\max_{1\leq j\leq q}|\hat{f}(\omega_{j})-f(\omega_{j})|=o_{\pr}(1/\log n).
\end{equation}
Note that under $H_{0}$, we have $Z_{n}=X_{n}$. For the briefness, we assume that $X_{n}$ satisfies $\ep |X_{n}|^{4+\gamma}<\infty$ for some
$\gamma>0$ and the geometric-moment contraction (GMC) condition $\theta_{n, 4+\gamma}=O(\rho^{n})$ for some $0<\rho<1$ holds. Many nonlinear time
series models, such as  GARCH models, generalized random coefficient autogressive models, nonlinear AR models, bilinear models, satisfy GMC; see
Section 5 in Shao and Wu \cite{shaox} for more details. By Lemma A.4 in Shao and Wu \cite{shaox}, we have
\begin{equation}
\max_{j,k\leq q}|\Cov(I_{n,X}(\omega_{k}),I_{n,X}(\omega_{j}))-f(\omega_{j})\delta_{j,k}|=O(1/n),
\end{equation}
where $\delta_{j,k}=I_{j=k}$, and it follows that
\begin{equation*}
q^{-1}\sum_{i=1}^{q}(I_{n,X}(\omega_{i})-\ep I_{n,X}(\omega_{i}))/f(\omega_{i})=O_{\pr}(1/\sqrt{n}).
\end{equation*}
Moreover, since $ I_{n,X}(\omega)=n^{-1}\sum_{k=-n+1}^{n-1}\sum_{t=1}^{n-|k|}X_{t}X_{t+|k|}\exp(-\text{i}k \omega), $
 we see that
$ \max_{\omega\in R}\Big{|}\frac{\ep I_{n,X}(\omega)}{2\pi f(\omega)}-1\Big{|}=O(1/n). $ This implies (\ref{re1}).

Now we choose the estimator
\begin{equation*}
\hat{f}(\omega)=\frac{1}{2\pi}\sum_{k=-B_{n}}^{B_{n}}\hat{r}(k)a(k/B_{n})\exp(-\text{i}k\omega),
\end{equation*}
where $\hat{r}(k)=n^{-1}\sum_{j=1}^{n-|k|}X_{j}X_{j+|k|}$, $|k|<n$, $a(\cdot)$ is an even, Lipschitz continuous function with support $[-1,1]$,
$a(0)=1$ and $a(x)-1=O(x^{2})$ as $x\rightarrow 0$, and $B_{n}$ is a sequence of positive integers with $B_{n}\rightarrow\infty$ and
$B_{n}/n\rightarrow 0$. Suppose now $B_{n}=O(n^{\eta})$, $0<\eta<\gamma/(4+\gamma)$, $0<\gamma<4$. Then Theorem 3.2 in Shao and Wu \cite{shaox}
gives
\begin{equation*}
\max_{\omega\in [0,\pi]}|\hat{f}(\omega)-\ep \hat{f}(\omega)|=O_{\pr}(\sqrt{B_{n}(\log n)/n}).
\end{equation*}
Moreover, simple calculations as in Woodroofe and Van Ness \cite{wood} imply $\max_{\omega\in [0,\pi]}|\ep
\hat{f}(\omega)-f(\omega)|=O(B^{-2}_{n})$. Hence (\ref{re2}) holds by letting $B_{n}\asymp n^{\eta}$, $0<\eta<\gamma/(4+\gamma)$. Finally, Theorem
2.2 together with (\ref{re1}) and (\ref{re2}) yields, under $H_{0}$, $g_{n}(Z)-\log q\Rightarrow G$, where $G$ has the standard Gumbel
distribution.

\section{Inequalities for Fourier transforms of stationary process}
 \setcounter{equation}{0}
In this section, we prove some inequalities for $X_{n}$ defined in (\ref{a0}). Suppose that $\ep X_{0}=0$ and $\ep X^{2}_{0}<\infty$. Note that
\begin{equation*}
X_{n}=\sum_{j\in Z}(\ep[X_{n}|\mathcal{F}_{-j,\infty}]-\ep[X_{n}|\mathcal{F}_{-j+1,\infty}])=:\sum_{j\in Z}\mathcal{P}_{j}(X_{n}).
\end{equation*}
By virtue of  H\"{o}lder's inequality, we have for $u\geq 0$,
\begin{align}\label{ab1.1}
|r(u)|=|\ep X_{0}X_{u}|=|\sum_{j\in Z}\ep \mathcal{P}_{j}(X_{0})\mathcal{P}_{j}(X_{u})|\leq \sum_{j=0}^{\infty}\theta_{j,2}\theta_{u+j,2},
\end{align}
and hence $ \sum_{u\geq n}|r(u)|\leq \Theta_{0,2}\Theta_{n,2}$.

Next, we approximate the Fourier transforms of $X_{n}$ by the sum of $m$-dependent random variables. Set
\begin{equation*}
X_{k}(m)=\ep[X_{k}|\varepsilon_{k-m},\cdots,\varepsilon_{k}],\quad k\in Z,~ m\geq 0.
\end{equation*}

\begin{lemma}\label{le1} Suppose that $\ep|X_{0}|^{p}<\infty$ for some $p\geq
2$ and $\Theta_{0,p}<\infty$. We have
\begin{equation*}
\sup_{\omega\in R}\ep\Big{|}\sum_{k=1}^{n}(X_{k}-X_{k}(m))\exp(\textbf{i}\omega k)\Big{|}^{p}\leq C_{p}n^{p/2}\Theta^{p}_{m,p},
\end{equation*}
where $C_{p}$ is a constant only depending on $p$.
\end{lemma}

{\bf Remark 3.1.} This lemma together with Proposition 1 in Wu \cite{wu1} would lead to the  maximal inequality: for $p>2$,
\begin{equation*}
\sup_{\omega\in R}\ep\max_{1\leq j\leq n}\Big{|}\sum_{k=1}^{j}(X_{k}-X_{k}(m))\exp(\text{i}\omega k)\Big{|}^{p}\leq C_{p}n^{p/2}\Theta^{p}_{m,p}.
\end{equation*}

\begin{proof} We decompose $X_{k}-X_{k}(m)$ as:
\begin{align*}
X_{k}-X_{k}(m)=\sum_{j=-k+m}^{\infty}(\ep[X_{k}|\mathcal {F}_{-j-1,k}] -\ep[X_{k}|\mathcal {F}_{-j,k}]) =:\sum_{j=-k+m}^{\infty}R_{k,j}.
\end{align*}
Therefore,
\begin{align*}
\sum_{k=1}^{n}\{X_{k}-X_{k}(m)\}\exp(\text{i}\omega k)=\sum_{j=-n+m}^{\infty}\sum_{k=1\vee(-j+m)}^{n}R_{k,j}\exp(\text{i}\omega k).
\end{align*}
For every fixed $n$ and $m$, $\{\sum_{k=1\vee(-j+m)}^{n}R_{k,j}\exp(\text{i}\omega k), j\geq -n+m\}$ is a sequence of martingale differences.
Hence by the Marcinkiewicz-Zygmund-Burkholder inequality,
\begin{align*}
&\ep\Big{|}\sum_{j=-n+m}^{\infty}\sum_{k=1\vee(-j+m)}^{n}R_{k,j}\exp(\text{i}\omega k)\Big{|}^{p}\leq
C_{p}\Big{(}\sum_{j=-n+m}^{\infty}\Big{(}\sum_{k=1\vee(-j+m)}^{n}\|R_{k,j}\|_{p}\Big{)}^{2}\Big{)}^{p/2}\\
&\leq C_{p}\Big{(}\sum_{j=-n+m}^{\infty}\Big{(}\sum_{k=1\vee(-j+m)}^{n}\theta_{j+1+k,p}\Big{)}^{2}\Big{)}^{p/2}\leq C_{p}n^{p/2}\Theta^{p}_{m,p}.
\end{align*}
This proves the lemma.
\end{proof}

Letting $m=0$ in Lemma \ref{le1} and noting that $X_{1}(0), X_{2}(0),\cdots$ are i.i.d. random variables, we obtain the following moment
inequalities.
\begin{lemma}\label{le2} Under the conditions of Lemma \ref{le1}, we
have, for $p\geq 2$,
\begin{equation*}
\ep\Big{|}\sum_{k=1}^{n}X_{k}\exp(\textbf{i}k\omega)\Big{|}^{p}\leq
Cn^{p/2}\mbox{~and~~}\ep\Big{|}\sum_{k=1}^{n}X_{k}(m)\exp(\textbf{i}k\omega)\Big{|}^{p}\leq Cn^{p/2},
\end{equation*}
where $C$ is a constant which does not depend on $\omega$ and $m$.
\end{lemma}

Define $S_{n,j,1}=\sum_{k=1}^{n}X_{k}\cos(k\omega_{j})$, $S_{n,j,2}=\sum_{k=1}^{n}X_{k}\sin(k\omega_{j})$, $1\leq j\leq q$.
\begin{lemma}\label{le3}
Suppose that $\ep X_{0}=0$, $\ep X^{2}_{0}<\infty$ and $\Theta_{0,2}<\infty$. Then
\\
(i).
\begin{equation*}
\max_{1\leq j\leq q}\Big{|}\frac{\ep S^{2}_{n,j,1}}{\pi nf(\omega_{j})}-1\Big{|}\leq Cn^{-1}\sum_{k=0}^{n}\Theta_{k,2}.
\end{equation*}
(ii).
\begin{equation*}
\max_{1\leq j\leq q}\Big{|}\frac{\ep S^{2}_{n,j,2}}{\pi nf(\omega_{j})}-1\Big{|}\leq Cn^{-1}\sum_{k=0}^{n}\Theta_{k,2}.
\end{equation*}
(iii).   $\max_{1\leq i, j\leq q}|\ep S_{n,i,1}S_{n,j,2}|\leq C\sum_{k=0}^{n}\Theta_{k,2}$ and $\max_{1\leq i\neq j\leq q}|\ep
S_{n,i,l}S_{n,j,l}|\leq C\sum_{k=0}^{n}\Theta_{k,2}$ for $l=1,2$.
\end{lemma}
\begin{proof} We only prove {\em (i)}, since the others can  be
obtained in an analogous way. We recall  the following propositions on the trigonometric functions:
\begin{itemize}
\item[(1)]{$\sum_{k=1}^{n}\cos(\omega_{j}k)\cos(\omega_{l}k)=\delta_{j,l}n/2$; (2)
~$\sum_{k=1}^{n}\sin(\omega_{j}k)\sin\omega_{l}k)=\delta_{j,l}n/2$;} \item[(3)]{$\sum_{k=1}^{n}\cos(\omega_{j}k)\sin(\omega_{l}k)=0$.}
\end{itemize}
By applying the above propositions, it is readily seen that
\begin{align*}
\frac{\ep S^{2}_{n,j,1}}{n}&=\frac{1}{2}\ep X^{2}_{1}+2n^{-1}\sum_{k=2}^{n}\sum_{i=1}^{k-1}\ep
X_{k}X_{i}\cos(k\omega_{j})\cos(i\omega_{j})\\
&=\frac{1}{2}\ep
X^{2}_{1}+2n^{-1}\sum_{k=1}^{n-1}r(k)\sum_{i=1}^{n-k}\cos(i\omega_{j})\cos((i+k)\omega_{j})\\
&=\frac{1}{2}\ep X^{2}_{1}+\sum_{k=1}^{n-1}r(k)\cos(k\omega_{j})\\
&\quad-2n^{-1}\sum_{k=1}^{n-1}r(k)\sum_{i=n-k+1}^{n}\cos(i\omega_{j})\cos((i+k)\omega_{j}),
\end{align*}
which, together with (\ref{ab1.1}) and the Abel lemma, implies
\begin{align*}
\Big{|}\frac{\ep S^{2}_{n,j,1}}{\pi nf(\omega_{j})}-1\Big{|}&\leq C\sum_{k=n}^{\infty}|r(k)|+Cn^{-1}\sum_{k=1}^{n-1}k|r(k)|\non &\leq
C\Theta_{n,2}+Cn^{-1}\sum_{j=0}^{\infty}\theta_{j,2}\sum_{k=1}^{n}k(\Theta_{k+j,2}-\Theta_{k+j+1,2})\\
&\leq Cn^{-1}\sum_{k=0}^{n}\Theta_{k,2}.
\end{align*}
 The proof of the lemma is complete.
\end{proof}

Let $m=[n^{\beta}]$ for some $0<\beta<1$ and $ J_{n,X}(\omega)=\Big{|}\sum_{k=1}^{n}\{X_{k}-X_{k}(m)\}\exp(\text{i}\omega k)\Big{|}. $
\begin{lemma}Suppose that $\ep X^{2}_{0}<\infty$ and
$\Theta_{n,2}=o(1/\log n)$. We have for any $0<\beta<1$,
\begin{equation*}
\max_{1\leq i\leq q}J_{n,X}(\omega_{i})=o_{\pr}(\sqrt{n/\log n}).
\end{equation*}
\end{lemma}
\begin{proof} Since $\Theta_{m,2}=o((\log n)^{-1})$, there exists a sequence  $\{\gamma_{n}\}$
 with $\gamma_{n}>0$ and $\gamma_{n}\rightarrow 0$ such that
$\Theta_{m,2}\leq\gamma_{n}(\log n)^{-1}$.  By the decomposition used in the proof of Lemma \ref{le1}, $
J_{n,X}(\omega)=|\sum_{j=-n+m}^{\infty}\sum_{k=1\vee(m-j)}^{n}R_{k,j}\exp(\text{i}k\omega)| $.
 Set
\begin{eqnarray*}
&&R_{j}(\omega)=\sum_{k=1\vee(m-j)}^{n}R_{k,j}\exp(\text{i}k\omega ),~\widetilde{R_{j}}(\omega)=R_{j}(\omega)I\Big{\{}|R_{j}(\omega)|\leq
\gamma_{n}\sqrt{\frac{n}{(\log n)^{3}}}\Big{\}},\\
&&\overline{R}_{j}(\omega)=\widetilde{R_{j}}(\omega)-\ep[\widetilde{R_{j}}(\omega)|\mathcal{F}_{-j,\infty}],
~\widehat{R}_{j}(\omega)=R_{j}(\omega)-\overline{R}_{j}(\omega).
\end{eqnarray*}
 Using the fact $\max_{\omega\in R}|R_{j}(\omega)|\leq
\sum_{k=1\vee(m-j)}^{n}|R_{k,j}|$, we see that for any $\delta>0$,
\begin{eqnarray*}
&&\pr\Big{(}\max_{\omega\in R}|\sum_{j=-n+m}^{\infty}\widehat{R}_{j}(\omega)|\geq \delta\sqrt{n/\log n}\Big{)}\leq C_{\delta}n^{-1/2}(\log
n)^{1/2}\sum_{j=-n+m}^{\infty}\ep\max_{\omega\in
R}|\widehat{R}_{j}(\omega)|\\
&&\leq 2C_{\delta}\frac{(\log n)^{2}\gamma^{-1}_{n}}{n}\sum_{j=-n+m}^{\infty}\Big{(}\sum_{k=1\vee (m-j)}^{n}\theta_{k+j+1,2}\Big{)}^{2}\leq
2C_{\delta}(\log n)^{2}\gamma^{-1}_{n}\Theta^{2}_{m,2}=o(1).
\end{eqnarray*}
Hence, in order to prove the lemma, it is sufficient to show that
\begin{eqnarray}\label{rresv}
\max_{1\leq i\leq q}|\sum_{j=-n+m}^{\infty}\overline{R}_{j}(\omega_{i})|=o_{\pr}(\sqrt{n/\log n}).
\end{eqnarray}
Setting the event $ A=\Big{\{}\max_{\omega\in R}\sum_{j=-n+m}^{\infty}\ep[|\overline{R}_{j}(\omega)|^{2}|\mathcal{F}_{-j,\infty}]\geq \gamma_{n}
n/(\log n)^{2}\Big{\}}$, we have
\begin{eqnarray*}
\pr(A) &\leq& C_{\delta}\frac{(\log n)^{2}\gamma^{-1}_{n}}{n}\sum_{j=-n+m}^{\infty}\ep\Big{(} \sum_{k=1\vee(m-j)}^{n}|R_{k,j}|\Big{)}^{2}\\&\leq&
C_{\delta}(\log n)^{2}\gamma^{-1}_{n}\Theta^{2}_{m,2}=o(1).
\end{eqnarray*}
Note that $\overline{R}_{j}(\omega)$, $j\geq -n+m$, are martingale differences. By applying  Freedman's inequality \cite{freed}, one concludes
that
\begin{eqnarray*}
\pr\Big{(}\max_{1\leq i\leq q}|\sum_{j=-n+m}^{\infty}\overline{R}_{j}(\omega_{i})|\geq \delta\sqrt{n/\log n}\Big{)}\leq
2n\exp\Big{(}-\frac{\delta^{2}\log n}{\gamma_{n}(8+8\delta)}\Big{)}+\pr(A)=o(1).
\end{eqnarray*}
This proves (\ref{rresv}).
\end{proof}

{\bf Remark 3.2.} Let $X_{n}=g((\varepsilon_{n-i})_{i\in Z})$ be a two-sided process.  For  $n\in Z$, denote $X^{*}_{n}$ by replacing
$\varepsilon_{0}$ with $\varepsilon^{*}_{0}$ in $X_{n}$. Define the physical dependence measure $\theta_{n,p}=\|X_{n}-X^{*}_{n}\|_{p}$ and
$\Theta_{n,p}=\sum_{|i|\geq n}\theta_{i,p}$. Also let $X_{k}(m)=\ep[X_{k}|\varepsilon_{k-m},\cdots,\varepsilon_{k+m}]$. Then Lemmas 3.1-3.4 still
hold for $X_{n}=g((\varepsilon_{n-i})_{i\in Z})$. This can be proved similarly by observing that
\begin{align}\label{b1}
X_{k}-X_{k}(m)&=\sum_{j=-k+m}^{\infty}(\ep[X_{k}|\mathcal {F}_{-j-1,\infty}] -\ep[X_{k}|\mathcal {F}_{-j,\infty}])\non
&\quad+\sum_{j=m+k}^{\infty}(\ep[X_{k}|\mathcal{F}_{k-m,j+1}]-\ep[X_{k}|\mathcal{F}_{k-m,j}])\non
&=:\sum_{j=-k+m}^{\infty}R^{(1)}_{k,j}+\sum_{j=m+k}^{\infty}R^{(2)}_{k,j},
\end{align}
$\|R^{(1)}_{k,j}\|_{p}\leq \theta_{k+j+1,p}$ and $\|R^{(2)}_{k,j}\|_{p}\leq \theta_{k-j-1,p}$. The details can be found in \cite{lin}.

\section{Proof of Theorem 2.1}
 \setcounter{equation}{0}
Let $h$ be a Lipschitz continuous function on $R$. Set
\begin{equation*}
\varepsilon^{'}_{i}=\varepsilon_{i}I\{|\varepsilon_{i}|\leq \gamma_{n}\sqrt{n/\log n}\}-\ep \varepsilon_{i}I\{|\varepsilon_{i}|\leq
\gamma_{n}\sqrt{n/\log n}\} , i\in Z,
\end{equation*}
 where $\gamma_{n}\rightarrow
0$. Put $Y^{'}_{k}=\sum_{i\in Z}a_{i}\varepsilon^{'}_{k-i}$, $X^{'}_{k}=h(Y^{'}_{k})-\ep h(Y^{'}_{k})$ for $1\leq k\leq n$. Since $\ep
\varepsilon^{2}_{0}I\{|\varepsilon_{0}|\geq n\}=o(1/\log n)$, we can choose $\gamma_{n}\rightarrow 0$  sufficiently slowly such that
\begin{equation*}
\sqrt{n\log n}\ep|\varepsilon_{0}|I\{|\varepsilon_{0}|\geq \gamma_{n}\sqrt{n/\log n}\}\rightarrow 0.
\end{equation*}
This together with the Lipschitz continuity of $h$ implies that
\begin{align*}
&\frac{\sqrt{\log n} \ep\max_{1\leq j\leq
q}|\sum_{k=1}^{n}(X_{k}-X^{'}_{k})\exp(\text{i}k\omega_{j})|}{\sqrt{n}}\\
&\leq C\sqrt{n\log n}\sum_{j\in Z}|a_{j}|\ep|\varepsilon_{0}|I\{|\varepsilon_{0}|\geq \gamma_{n}\sqrt{n/\log n}\}\rightarrow 0.
\end{align*}
In addition, note that for $1\leq j\leq q$,
\begin{align*}
|I_{n,X}(\omega_{j})-I_{n,X^{'}}(\omega_{j})|&\leq \sqrt{M_{n}(X^{'})}\max_{1\leq j\leq
q}|\sum_{k=1}^{n}(X_{k}-X^{'}_{k})\exp(\text{i}k\omega_{j})|/\sqrt{n}\\
&\quad+\max_{1\leq j\leq q}|\sum_{k=1}^{n}(X_{k}-X^{'}_{k})\exp(\text{i}k\omega_{j})|^{2}/n.
\end{align*}
Then, in order to prove Theorem 2.1, we only need to show that
\begin{equation*}
I_{n,q}(X^{'})-\log q\Rightarrow G.
\end{equation*}
 Recall that $m=[n^{\beta}]$ for some $0<\beta<1$. Let
\begin{equation*}
X^{'}_{k}(m)=\ep[X^{'}_{k}|\varepsilon_{k-m},\cdots, \varepsilon_{k+m}],\quad 1\leq k\leq n,
\end{equation*}
and
\begin{equation*}
\widetilde{J}_{n,X}(\omega)=\Big{|}\sum_{k=1}^{n}(X^{'}_{k}-X^{'}_{k}(m))\exp(\text{i}\omega k)\Big{|}.
\end{equation*}
By Lemma 3.4 and Remark 3.2, it is readily seen that
\begin{eqnarray}\label{revi}
\max_{1\leq i\leq q}\widetilde{J}_{n,X}(\omega_{i})=o_{\pr}(\sqrt{n/\log n}).
\end{eqnarray}
We define the periodogram $ I_{n,X^{'}(m)}(\omega)=n^{-1}\Big{|}\sum_{k=1}^{n}X^{'}_{k}(m)\exp(\text{i}k\omega)\Big{|}^{2}, $ and let
$I_{n,1}(X^{'}(m))\leq\cdots\leq I_{n,q}(X^{'}(m))$ be
 the order statistics of $I_{n,X^{'}(m)}(\omega_{j})/(2\pi f(\omega_{j}))$, $1\leq j\leq q$. In view of (\ref{revi}),
 it is sufficient to prove  that
\begin{equation}\label{a6}
I_{n,q}(X^{'}(m))-\log q\Rightarrow G.
\end{equation}

For $0<\beta<\alpha<1/10$, let us split the interval $[1, n]$
 into
\begin{align*}
&H_{j}=[(j-1)(n^{\alpha}+2n^{\beta})+1,
(j-1)(n^{\alpha}+2n^{\beta})+n^{\alpha}],\\
&I_{j}=[(j-1)(n^{\alpha}+2n^{\beta})+n^{\alpha}+1,
j(n^{\alpha}+2n^{\beta})],\\
&1\leq j\leq m_{n}-1,\quad m_{n}-1=[n/(n^{\alpha}+2n^{\beta}))]\sim
n^{1-\alpha},\\
&H_{m_{n}}=[(m_{n}-1)(n^{\alpha}+2n^{\beta})+1, n].
\end{align*}
Here and below the  notation $n^{\alpha}$ is used to denote $[n^{\alpha}]$ for briefness. Put $v_{j}(\omega)=\sum_{k\in
I_{j}}X^{'}_{k}(m)\exp(\text{i}k\omega)$, $1\leq j\leq m_{n}-1$. Then $v_{j}(\omega)$, $1\leq j\leq m_{n}-1$, are independent and can be neglected
by observing the following lemma.
\begin{lemma}Under (\ref{a3}), we have $\max_{1\leq l\leq q}|\sum_{j=1}^{m_{n}-1}v_{j}(\omega_{l})|=o_{\pr}(\sqrt{n/\log
n})$.
\end{lemma}
\begin{proof}
First,   Corollary 1.6 of Nagaev \cite{nagaev}, which is a Fuk-Nagaev-type inequality, shows that for any large $Q$,
\begin{align*}
&\sum_{l=1}^{q}\pr\Big{(}|\sum_{j=1}^{m_{n}-1}v_{j}(\omega_{l})|\geq\delta\sqrt{n/\log
n}\Big{)}\\
&\leq C_{Q,\delta}\sum_{l=1}^{q}\Big{(}\frac{\sum_{j=1}^{m_{n}-1}\ep v^{2}_{j}(\omega_{l})}{n/\log
n}\Big{)}^{Q}+C_{Q}\sum_{l=1}^{q}\sum_{j=1}^{m_{n}-1}\pr\Big{(}|v_{j}(\omega_{l})|\geq C_{Q}\delta\sqrt{n/\log n}\Big{)}.
\end{align*}
By Lemma \ref{le2} and Remark 3.2, $\sum_{j=1}^{m_{n}-1}\ep v^{2}_{j}(\omega_{l})\leq Cn^{1-\alpha+\beta}$. So the first term above tends to zero.
To complete the proof of Lemma 4.1, we shall show the second term also tends to zero. In fact, using the fact $|h(x)|\leq C(|x|+1)$, we can get
\begin{align}\label{b6}
|v_{j}(\omega_{l})| &\leq C\Big{|}\sum_{k\in I_{j}}\sum_{i=-m}^{m}|a_{i}|(|\varepsilon^{'}_{k-i}|-\ep|\varepsilon^{'}_{k-i}|)\Big{|}+C|I_{j}|\non
&=_{d}C\Big{|}\sum_{k\in I_{1}}\sum_{i=-m}^{m}|a_{i}|(|\varepsilon^{'}_{k-i}|-\ep|\varepsilon^{'}_{k-i}|)\Big{|}+C|I_{1}|\non
&=C\Big{|}\sum_{t=-m}^{3m}\sum_{k=1\vee (t-m)}^{(m+t)\wedge (2m)}|a_{k-t}|(|\varepsilon^{'}_{t}|-\ep|\varepsilon^{'}_{t}|)\Big{|}+C|I_{1}|,
\end{align}
where $X=_{d}Y$ means $X$ and $Y$ have the same distribution.  Hence
\begin{align}\label{b4}
&\sum_{l=1}^{q}\sum_{j=1}^{m_{n}-1}\pr\Big{(}|v_{j}(\omega_{l})|\geq C_{Q}\delta\sqrt{n/\log n}\Big{)}\non &\leq
\sum_{l=1}^{q}\sum_{j=1}^{m_{n}-1}\pr\Big{(}\Big{|}\sum_{t=-m}^{3m}\sum_{k=1\vee (t-m)}^{(m+t)\wedge
(2m)}|a_{k-t}|(|\varepsilon^{'}_{t}|-\ep|\varepsilon^{'}_{t}|)\Big{|}\geq C_{Q}\delta\sqrt{n/\log n}\Big{)}\non
 &\leq C\sum_{l=1}^{q}\sum_{j=1}^{m_{n}-1}\Big{(}\frac{m}{n/\log
n}\Big{)}^{Q}\rightarrow 0,
\end{align}
where the last inequality follows from  the Fuk-Nagaev inequality, by noting that $|\varepsilon^{'}_{t}|\leq \gamma_{n}\sqrt{n/\log n}$. The
desired conclusion  is established.
\end{proof}
We now deal with the sum of large blocks. Let
\begin{align*}
 &u_{j}(\omega)=\sum_{k\in H_{j}}X^{'}_{k}(m)\exp(\text{i}k\omega),~u^{'}_{j}(\omega)=u_{j}(\omega)I\{|u_{j}(\omega)|\leq\gamma^{1/2}_{n}\sqrt{n/\log
n}\},\\
&\overline{u}_{j}(\omega)=u^{'}_{j}(\omega)-\ep u^{'}_{j}(\omega), ~ 1\leq j\leq m_{n}.
\end{align*}
Noting that $|u_{j}(\omega)|\leq \sum_{k\in H_{j}}|X^{'}_{k}(m)|=:\xi_{j}$, $m_{n}\sim n^{1-\alpha}$ and using similar arguments to those employed
in (\ref{b6}) and (\ref{b4}), it is readily seen that for any large $Q$,
\begin{align}\label{ab3}
&\frac{\sqrt{\log n}\sum_{j=1}^{m_{n}}\ep \xi_{j}I\{\xi_{j}\geq\gamma^{1/2}_{n}\sqrt{n/\log n}\} }{\sqrt{n}}\non &\quad\leq C\sqrt{\log
n}n^{1/2-\alpha}\sum_{k=n}^{\infty}\frac{1}{\sqrt{k\log k}}\pr\Big{(}\xi_{1}\geq\gamma^{1/2}_{n}\sqrt{k/\log k}\Big{)}\non
&\quad\quad+Cn^{1-\alpha}\pr\Big{(}\xi_{1}\geq \gamma^{1/2}_{n}\sqrt{n/\log n}\Big{)}\non
 &\quad\leq C\sqrt{\log
n}n^{1/2-\alpha}\sum_{k=n}^{\infty}\frac{1}{\sqrt{k\log k}}\Big{(}\frac{n^{\alpha}}{\gamma_{n}k/\log k}\Big{)}^{Q}\non
&\quad\quad+Cn^{1-\alpha}(\gamma^{-1}_{n}n^{\alpha-1}\log n)^{Q}\non &\quad=o(1),
\end{align}
which implies $ \max_{1\leq l\leq q}|\sum_{j=1}^{m_{n}}\Big{(}u_{j}(\omega_{l})-\overline{u}_{j}(\omega_{l})\Big{)}|=o_{\pr}(\sqrt{n/\log n}). $
Combining  this and Lemma 4.1 yields that we only need to show
\begin{equation}\label{b5}
I_{n,q}(\overline{X})-\log q\Rightarrow G,
\end{equation}
where $I_{n,q}(\overline{X})$ denotes the maximum of
$$|\sum_{k=1}^{m_{n}}\overline{u}_{k}(\omega_{l})|^{2}/(2\pi
nf(\omega_{l})), 1\leq l\leq q.$$

In order to prove (\ref{b5}), we need the following moderate deviation result, whose proof is
 based on Gaussian approximation technique due to Einmahl \cite{einmahl}, Corollary 1(b), page
 31 and Remark on page 32. The detailed proof is given in \cite{lin}.

 \begin{lemma} Let $\xi_{n,1},\cdots, \xi_{n,k_{n}}$ be independent random vectors with mean zero
 and values in $\mathbb{R}^{2d}$, and $S_{n}=\sum_{i=1}^{k_{n}}\xi_{n,i}$. Assume that $|\xi_{n,k}|\leq
 c_{n}B^{1/2}_{n}$, $1\leq k\leq k_{n}$, for some $c_{n}\rightarrow 0$,
 $B_{n}\rightarrow\infty$ and
 \begin{equation*}
 \Big{|}B^{-1}_{n}\Cov(\xi_{n,1}+\cdots+\xi_{n,k_{n}})-I_{2d}\Big{|}=O(c^{2}_{n}),
 \end{equation*}
 where $I_{2d}$ is a $2d\times 2d$ identity matrix. Suppose that
$
 \beta_{n}:=B^{-3/2}_{n}\sum_{k=1}^{k_{n}}\ep|\xi_{n,k}|^{3}\rightarrow
 0.
 $
 Then
 \begin{align*}
& |\pr(|S_{n}|_{2d}\geq x)-\pr(|N|_{2d}\geq
x/B^{1/2}_{n})|\\
&\leq o(\pr(|N|_{2d}\geq x/B^{1/2}_{n}))+C\Big{(}\exp\Big{(}-\frac{\delta^{2}_{n}\min(c^{-2}_{n}, \beta^{-2/3}_{n})}{16d}\Big{)}
+\exp\Big{(}\frac{Cc^{2}_{n}}{\beta^{2}_{n}\log \beta_{n}}\Big{)}\Big{)},
 \end{align*}
uniformly for $x\in[B^{1/2}_{n}, \delta_{n}\min(c^{-1}_{n}, \beta^{-1/3}_{n})B^{1/2}_{n}]$, with any $\delta_{n}\rightarrow 0$ and
$\delta_{n}\min(c^{-1}_{n}, \beta^{-1/3}_{n})\rightarrow \infty$. $N$ is a centered normal random vector with covariance matrix $I_{2d}$.
$|\cdot|_{2d}$ is defined by $|z|_{2d}=\min\{(x^{2}_{i}+y^{2}_{i})^{1/2}: 1\leq i\leq
  d\}$, $z=(x_{1}, y_{1},\cdots, x_{d},y_{d})$.
 \end{lemma}

We begin the proof of (\ref{b5}) by checking the conditions in Lemma 4.2. We  define the following notations:
$\overline{u}_{k}(\omega_{l})/f^{1/2}(\omega_{l})=:\overline{u}_{k,l}(1)+\text{i}\overline{u}_{k,l}(2)$,
\begin{equation}\label{fh}
Z_{k}=(\overline{u}_{k,i_{1}}(1), \overline{u}_{k,i_{1}}(2),\cdots, \overline{u}_{k,i_{d}}(1), \overline{u}_{k,i_{d}}(2)),~1\leq
i_{1}<\cdots<i_{d}\leq q
\end{equation}
 and $U_{n}=\sum_{k=1}^{m_{n}}Z_{k}$. Then it is easy to see that $Z_{1},\cdots, Z_{m_{n}}$ are independent.
 \begin{lemma}Under the conditions of Theorem 2.1, we have
 \begin{equation*}
\Big{|}\Cov(U_{n})/(n\pi)-I_{2d}\Big{|}=o(1/\log n)
 \end{equation*}
  uniformly for $1\leq
i_{1}<\cdots<i_{d}\leq q$.
 \end{lemma}\vspace{-5mm}
 \begin{proof}
Let $B_{n,i}=\sum_{k=1}^{m_{n}}\ep (\overline{u}_{k,i}(1))^{2}$. Similar arguments to those in (\ref{ab3}) together with some elementary
calculations give that $\max_{1\leq l\leq q}\ep |u_{j}(\omega_{l})-\overline{u}_{j}(\omega_{l})|^{2}=O(n^{-Q})$ for any large $Q$. This yields
that, for any large $Q$,
\begin{align}\label{aq1}
&\Big{|}B_{n,i}-\sum_{j=1}^{m_{n}}\ep\Big{(}\sum_{k\in H_{j}}X^{'}_{k}(m)\cos(k\omega_{i})\Big{)}^{2}\Big{|}\non &\quad\leq
C\sum_{j=1}^{m_{n}}|H_{j}|^{1/2}(\ep |u_{j}(\omega_{i})-\overline{u}_{j}(\omega_{i})|^{2})^{1/2}+\sum_{j=1}^{m_{n}}\ep
|u_{j}(\omega_{i})-\overline{u}_{j}(\omega_{i})|^{2}\non &\quad\leq Cn^{-Q}.
\end{align}
 Moreover, it follows from Lemmas \ref{le2}
and \ref{le1} and Remark 3.2 that
\begin{eqnarray}\label{aq2}
&&\Big{|}\ep\Big{(}\sum_{k=1}^{n}X^{'}_{k}(m)\cos(k\omega_{i})\Big{)}^{2}-\sum_{j=1}^{m_{n}}\ep\Big{(}\sum_{k\in
H_{j}}X^{'}_{k}(m)\cos(k\omega_{i})\Big{)}^{2}\Big{|}\leq Cn^{1-(\alpha-\beta)/2},\non
&&\quad\quad\Big{|}\ep\Big{(}\sum_{k=1}^{n}X^{'}_{k}(m)\cos(k\omega_{i})\Big{)}^{2}-
\ep\Big{(}\sum_{k=1}^{n}X^{'}_{k}\cos(k\omega_{i})\Big{)}^{2}\Big{|}=o(n/\log n).
\end{eqnarray}
In the  case $h(x)\equiv x$, we have $ \sum_{k=1}^{n}X^{'}_{k}\cos(k\omega_{i})=\sum_{t=-\infty}^{\infty}\sum_{k=1}^{n}a_{k+t}\cos(k\omega_{i})
\varepsilon^{'}_{-t}. $ Hence, condition (\ref{a4}) ensures that
\begin{equation}\label{ab4}
\Big{|} \ep\Big{(}\sum_{k=1}^{n}X^{'}_{k}\cos(k\omega_{i})\Big{)}^{2}- \ep\Big{(}\sum_{k=1}^{n}X_{k}\cos(k\omega_{i})\Big{)}^{2}\Big{|}=o(n/\log
n).
\end{equation}
Suppose now that $h$ is  Lipschitz continuous. We write $\zeta_{k}=|\varepsilon_{k}|I\{|\varepsilon_{k}|\geq \gamma_{n}\sqrt{n/\log n}\}$.  Then,
since $ |X_{k}-X^{'}_{k}|\leq C\sum_{j\in Z}|a_{j}|(\zeta_{k-j}+\ep\zeta_{k-j})$, we have from $\ep\varepsilon^{2}_{0}I\{|\varepsilon_{0}|\geq
n\}=o(1/(\log n)^{2})$ and the fact $\gamma_{n}\rightarrow 0$ sufficiently slowly that
\begin{align*}
&\ep\Big{(}\sum_{k=1}^{n}(X_{k}-X^{'}_{k})\cos(k\omega_{i})\Big{)}^{2}\\
&\leq C\ep\Big{(}\sum_{k=1}^{n}\sum_{j\in Z}|a_{j}|(\zeta_{k-j}-\ep \zeta_{k-j})\Big{)}^{2}+C\Big{(}\sum_{k=1}^{n}\sum_{j\in
Z}|a_{j}|\ep \zeta_{k-j}\Big{)}^{2}\\
&\leq Cn\ep \zeta^{2}_{0}+Cn^{2}(\ep \zeta_{0})^{2}=o(n/(\log n)^{2}),
\end{align*}
which implies (\ref{ab4}) by virtue of  Lemma 3.2 and the inequality $|\ep X^{2}-\ep Y^{2}|\leq \|X-Y\|_{2}\|X+Y\|_{2}$ for any random variables
$X$ and $Y$. From Lemma \ref{le3}, Remark 3.2 and (\ref{aq1})-(\ref{ab4}), we have $|B_{n,i}/(n\pi)-1|=o(1/\log n)$ uniformly for $1\leq i\leq q$.

In the following, we show that the off-diagonal elements in $\Cov(U_{n})$ are  $o(n/\log n)$.  We only deal with $B_{n,i,j}:=\ep
\{\sum_{k=1}^{m_{n}} \overline{u}_{k,i}(1)\sum_{k=1}^{m_{n}} \overline{u}_{k,j}(1)\}$, $i\neq j$, since the other elements can be estimated
similarly. As in (\ref{aq1}) and (\ref{aq2}), we have
\begin{align*}
&\Big{|}B_{n,i,j}-(f(\omega_{i})f(\omega_{j}))^{-\frac{1}{2}}\ep\Big{(}\sum_{k=1}^{n}X^{'}_{k}(m)\cos(k\omega_{i})\sum_{k=1}^{n}X^{'}_{k}(m)\cos(k\omega_{j})\Big{)}\Big{|}\\
&\leq C\Big{|}\ep \Big{[}\Big{\{}\sum_{k=1}^{m_{n}}
\overline{u}_{k,i}(1)-(f(\omega_{i}))^{-\frac{1}{2}}\sum_{k=1}^{n}X^{'}_{k}(m)\cos(k\omega_{i})\Big{\}}\sum_{k=1}^{m_{n}}
\overline{u}_{k,j}(1)\Big{]}\Big{|}\\
&\quad+C|f(\omega_{i})|^{-\frac{1}{2}}\Big{|}\ep \Big{[}\sum_{k=1}^{n}X^{'}_{k}(m)\cos(k\omega_{i})\Big{\{}\sum_{k=1}^{m_{n}}
\overline{u}_{k,j}(1)-(f(\omega_{j}))^{-\frac{1}{2}}\sum_{k=1}^{n}X^{'}_{k}(m)\cos(k\omega_{j})\Big{\}}\Big{]}\Big{|}\\
&\leq Cn^{1-(\alpha-\beta)/2}.
\end{align*}
Moreover, by virtue of Lemmas 3.1-3.3 and Remark 3.2, we have
\begin{align*}
&\ep\Big{(}\sum_{k=1}^{n}X^{'}_{k}(m)\cos(k\omega_{i})\sum_{k=1}^{n}X^{'}_{k}(m)\cos(k\omega_{j})\Big{)}=o(n/\log n).
\end{align*}
Hence $B_{n,i,j}=o(n/\log n)$, $i\neq j$. This proves the lemma.
\end{proof}
\begin{lemma}\label{lle}Under the conditions of Theorem 2.1, we have
uniformly for $1\leq i_{1}<\cdots<i_{d}\leq q$ that
\begin{equation*}
\overline{\beta}_{n}:=n^{-3/2}\sum_{j=1}^{m_{n}}\ep |Z_{j}|^{3}=o(1/(\log n)^{3/2}).
\end{equation*}
\end{lemma}
\begin{proof}
By the arguments in (\ref{b6}), the Fuk-Nagaev inequality and the fact $\alpha<1/10$ and  $\gamma_{n}\rightarrow 0$ sufficiently slowly,
\begin{align*}
\sum_{j=1}^{m_{n}}\ep |\overline{u}_{j}(\omega_{i})|^{3} &\leq \sum_{j=1}^{m_{n}}\sum_{k=1}^{n}\Big{(}\frac{k}{\log
k}\Big{)}^{3/2}\pr\Big{(}\gamma^{1/2}_{n}\sqrt{\frac{k}{\log k}}< |u_{j}(\omega_{i})|\leq \gamma^{1/2}_{n}\sqrt{\frac{k+1}{\log
(k+1)}}\Big{)}\\
&\leq Cn^{1+5\alpha}+C\sum_{j=1}^{m_{n}}\sum_{k=n^{4\alpha}}^{n}\frac{k^{1/2}}{(\log k)^{3/2}}\pr\Big{(} |u_{j}(\omega_{i})|\geq
\gamma^{1/2}_{n}\sqrt{\frac{k}{\log k}}\Big{)}\\
&\quad+ C\sum_{j=1}^{m_{n}}\frac{n^{6\alpha}}{(\log n)^{3/2}} \pr\Big{(} |u_{j}(\omega_{i})|\geq \gamma^{1/2}_{n}\sqrt{\frac{n^{4\alpha}}{\log
n^{4\alpha}}}\Big{)}
\\
&\leq Cn^{1+5\alpha}+C\sum_{j=1}^{m_{n}}\sum_{k=n^{4\alpha}}^{n} \frac{k^{1/2}}{(\log k)^{3/2}}\Big{(}\frac{n^{\alpha}}{\gamma_{n}k/\log
k}\Big{)}^{Q}\\
&\quad+C\sum_{j=1}^{m_{n}}\sum_{k=n^{4\alpha}}^{n} \frac{k^{1/2}n^{\alpha}}{(\log k)^{3/2}}\pr\Big{(}|\varepsilon_{0}|\geq
C\gamma^{1/2}_{n}\sqrt{\frac{k}{\log k}}\Big{)}\\
&\quad+C\sum_{j=1}^{m_{n}}\frac{n^{7\alpha}}{(\log n)^{3/2}}\pr\Big{(}|\varepsilon_{0}|\geq
C\gamma^{1/2}_{n}\sqrt{\frac{n^{4\alpha}}{\log n^{4\alpha}}}\Big{)}\\
 &=o((n/\log n)^{3/2}),\mbox{~uniformly for $1\leq i\leq q$.}
\end{align*}
The desired result now follows.
\end{proof}

 By Lemmas 4.3 and \ref{lle}, we may write $\overline{\beta}_{n}=\nu^{3/2}_{n}(\log
 n)^{-3/2}$ and $\Big{|}\Cov(U_{n})/(n\pi)-I_{2d}\Big{|}=\gamma_{n,1}(\log
 n)^{-1}$,
 where $\nu_{n}\rightarrow0,\gamma_{n,1}\rightarrow 0$.  Let us take $c_{n}=\{(4d\gamma_{n}(\pi f^{*})^{-1})^{1/2}\vee\gamma^{1/2}_{n,1}\}(\log
n)^{-1/2}=:\gamma^{1/2}_{n,2}(\log n)^{-1/2}$ and $\delta_{n}=\max\{\gamma^{1/4}_{n,2},\nu^{1/4}_{n}\}$ in Lemma  4.2. Note that
$\gamma_{n,2}\rightarrow 0$ sufficiently slowly. Then, simple calculations show that
\begin{equation*}
\exp\Big{(}-\frac{\delta^{2}_{n}\min(c^{-2}_{n}, \overline{\beta}^{-2/3}_{n})}{16d}\Big{)}\leq Cn^{-4d},\quad
\exp\Big{(}\frac{Cc^{2}_{n}}{\overline{\beta}^{2}_{n}\log \overline{\beta}_{n}}\Big{)}\leq Cn^{-4d}.
\end{equation*}
By virtue of Lemma 4.2, it holds that for any fixed $x\in R$,
\begin{align}\label{u1}
&\pr\Big{(}(2n\pi)^{-1/2}|U_{n}|_{2d}\geq x+\log q\Big{)}\non &\quad\quad=\pr(|N|_{2d}\geq \sqrt{2}(x+\log q))(1+o(1))\non
&\quad\quad=q^{-d}\exp(-dx)(1+o(1)),
\end{align}
uniformly for $1\leq i_{1}<\cdots<i_{d}\leq q$. We write $V_{j}:=|\sum_{k=1}^{m_{n}}\overline{u}_{k}(\omega_{j})|^{2}/(2\pi nf(\omega_{j}))$,
$1\leq j\leq q$, and
\begin{align*}
A&:=\{I_{n,q}(\overline{X})\geq x+\log q\}=\bigcup_{j=1}^{q}\{V_{i}\geq x+\log q\} =:\bigcup_{j=1}^{q}A_{j}.
\end{align*}
 By
the Bonferroni inequality, we have for any fixed $k$ satisfying $1\leq k\leq q$,
\begin{equation*}
\sum_{t=1}^{2k}(-1)^{t-1}E_{t}\leq\pr(A)\leq \sum_{t=1}^{2k-1}(-1)^{t-1}E_{t},
\end{equation*}
where $ E_{t}=\sum_{1\leq i_{1}<\cdots<i_{t}\leq q}\pr(A_{i_{1}}\cap\cdots\cap A_{i_{t}}).$ In view of (\ref{u1}), it follows that
$\lim_{n\rightarrow\infty}E_{t}=e^{-tx}/t!$. Since $\sum_{t=1}^{k}(-1)^{t-1}e^{-tx}/t!\rightarrow 1-e^{-e^{-x}}$ as $k\rightarrow\infty$, the
proof of Theorem 2.1 is complete.

\section{Proof of Theorem 2.2}
 \setcounter{equation}{0}
Recall that $m=[n^{\beta}]$ and $\beta$ is sufficiently small. Let $S_{n,m}(\omega)=\sum_{k=1}^{n}X_{k}(m)\exp(\text{i}\omega k)$ and
 $I_{n,1}(m)\leq\cdots\leq I_{n,q}(m)$ be
 the order statistics of
 $|S_{n,m}(\omega_{j})|^{2}/(2\pi nf(\omega_{j}))$, $1\leq j\leq q$.
By Lemma 3.4 we only need to prove that
\begin{equation}\label{ab6}
I_{n,q}(m)-\log q\Rightarrow G.
\end{equation}
We use the same notations and blocking method as in the proof of Theorem 2.1 (replacing $X^{'}_{k}(m)$ with $X_{k}(m)$). For example,
$v_{j}(\omega)=\sum_{k\in I_{j}}X_{k}(m)\exp(ik\omega)$. As in Lemma 4.1, we claim that
\begin{equation}\label{ab5}
\max_{1\leq j\leq q}\Big{|}\sum_{k=1}^{m_{n}-1}v_{k}(\omega_{j})\Big{|}=o_{\pr}(\sqrt{n/\log n}).
\end{equation}
We come to prove it. Recall that $s>2$ and $\beta<\alpha$. Then we can choose $\alpha, \beta$ sufficiently small and $\tau$ sufficiently close to
$1/2$ such that
\begin{eqnarray}\label{revis}
(s-1)^{-1}(1-\alpha+\alpha s-1/2)<\tau<1/2.
\end{eqnarray}
We define $\overline{v}_{k}(\omega_{j})=v^{'}_{k}(\omega_{j})-\ep v^{'}_{k}(\omega_{j})$, where $
v^{'}_{k}(\omega_{j})=v_{k}(\omega_{j})I\{|v_{k}(\omega_{j})|\leq n^{\tau}\},$  $1\leq j\leq q$, $1\leq k\leq m_{n}-1$. So
\begin{equation*}
\max_{1\leq j\leq q}|\sum_{k=1}^{m_{n}-1}v_{k}(\omega_{j})|\leq \max_{1\leq j\leq q}|\sum_{k=1}^{m_{n}-1}\overline{v}_{k}(\omega_{j})|+
\max_{1\leq j\leq q}|\sum_{k=1}^{m_{n}-1}(v_{k}(\omega_{j})-\overline{v}_{k}(\omega_{j}))|.
\end{equation*}
By the Fuk-Nagaev inequality and Lemma \ref{le2}, we have for any large $Q$,
\begin{align}\label{b9}
\pr\Big{(}\max_{1\leq j\leq q}|\sum_{k=1}^{m_{n}-1}\overline{v}_{k}(\omega_{j})|\geq \delta\sqrt{\frac{n}{\log n}}\Big{)}\leq
Cn\Big{(}\frac{n^{1-\alpha+\beta}}{n/\log n}\Big{)}^{Q} \rightarrow 0.
\end{align}
Also, using (\ref{revis}), the condition $\ep|X_{0}|^{s}<\infty$ and $|v_{k}(\omega)|\leq \sum_{j\in I_{k}}|X_{j}(m)|$, we can get
\begin{eqnarray}\label{reviss}
&&\ep \frac{\max_{1\leq j\leq q}|\sum_{k=1}^{m_{n}-1}(v_{k}(\omega_{j})-\overline{v}_{k}(\omega_{j}))|}{\sqrt{n/\log n}}\non &&
\leq\frac{2n^{1-\alpha}\ep\Big{[} \sum_{k=1}^{n^{\beta}}|X_{k}(m)|I\{\sum_{k=1}^{n^{\beta}}|X_{k}(m)|\geq n^{\tau}\}\Big{]}}{\sqrt{n/\log n}}\non
&&\leq Cn^{1-\alpha+\beta s-\tau(s-1)-1/2}(\log n)^{1/2}=o(1).
\end{eqnarray}
This together with (\ref{b9}) implies (\ref{ab5}).

Set
\begin{align*}
&u^{'}_{k}(\omega_{j})=u_{k}(\omega_{j})I\{|u_{k}(\omega_{j})|\leq
n^{\tau}\},\\
& \overline{u}_{k}(\omega_{j})=u^{'}_{k}(\omega_{j})-\ep u^{'}_{k}(\omega_{j}),~1\leq j\leq q, ~ 1\leq k\leq m_{n}.
\end{align*}
By the similar arguments as (\ref{reviss}), using (\ref{revis}), we can show that
\begin{equation*}
\max_{1\leq j\leq q}\Big{|}\sum_{k=1}^{m_{n}}(u_{k}(\omega_{j})-\overline{u}_{k}(\omega_{j}))\Big{|}=o_{\pr}(\sqrt{n/\log n}).
\end{equation*}
So in order to get (\ref{ab6}), similarly to (\ref{b5}), it is sufficient to prove
\begin{equation}\label{rrre}
I_{n,q}(\overline{X})-\log q\Rightarrow G.
\end{equation}
In fact, (\ref{rrre}) follows from Lemmas 5.1 and 5.2 and similar arguments to those employed in the proof of Theorem 2.1.
\begin{lemma}
Under the conditions of Theorem 2.2, we have
 \begin{equation*}
\Big{|}\Cov(U_{n})/(n\pi)-I_{2d}\Big{|}=o(1/\log n).
 \end{equation*}
\end{lemma}
\begin{proof}
The same arguments as those of Lemma 4.3 give that
\begin{equation*}
|B_{n,i}-\ep\Big{(}\sum_{k=1}^{n}X_{k}\cos(k\omega_{i})\Big{)}^{2}/(\pi f(\omega_{i}))|=o(n/\log n).
\end{equation*}
 The lemma then follows from Lemma \ref{le3}.
\end{proof}
\begin{lemma}Under the conditions of Theorem 2.2, we have
\begin{equation*}
\overline{\beta}_{n}=n^{-3/2}\sum_{j=1}^{m_{n}}\ep |Z_{j}|^{3}=O(n^{t-1/2}),
\end{equation*}
where $t=\max\{(3-s)\tau+\alpha(s-2)/2, \alpha/2\}<\tau<1/2$.
\end{lemma}
\begin{proof} Suppose that $2<s<3$. Then by virtue of Lemma \ref{le2}, we have
\begin{align*}
\overline{\beta}_{n}\leq Cn^{-3/2+(3-s)\tau}\sum_{j=1}^{m_{n}}\ep |Z_{j}|^{s}\leq Cn^{-3/2+(3-s)\tau}\sum_{j=1}^{m_{n}}|H_{j}|^{s/2}\leq
Cn^{t-1/2}.
\end{align*}
The case of $s\geq 3$ can be similarly proved.
\end{proof}

{\bf{Acknowledgements }} The authors would like to thank an Associate Editor and the referee for many valuable comments.

\title{\bf Supplementary Material for "On maxima of periodograms  of stationary processes"}
\date{ }

\maketitle
\begin{center}
\vskip -2cm { {\sc Zhengyan Lin\footnote{Email: zlin@zju.edu.cn}}, {\sc Weidong Liu\footnote {Email: liuweidong99@gmail.com }}
{\small\centerline{Department of Mathematics, Zhejiang University, Hangzhou 310027, China}} }
\end{center}

\bigskip

\bigskip
{\bf Abstract.} This note is the supplementary material for "On maxima of periodograms  of stationary processes" by Lin and Liu (2009).

\section{Proof of Lemma 3.3 for two-sided process} \setcounter{equation}{0} \noindent
Let $\{\varepsilon_{n}; n\in Z\}$ be independent and identically distributed (i.i.d.) random variables and $g$ be a measurable function such that
 \begin{equation}\label{a0}
 X_{n}=g((\varepsilon_{n-i})_{i\in Z})
 \end{equation}
is a well-defined random variable. In Lin and Liu (2007), we proved some inequalities for Fourier transforms of one-sided causal processes; see
Lemmas 3.1-3.4 in Lin and Liu (2007). The proof can be similarly extended to the two-sided stationary process $X_{n}$. The proofs of Lemmas
3.1,3.2 and 3.4 for $X_{n}$ in (\ref{a0}) are exactly the same as those given in Lin and Liu (2009), and so we do not repeat them. We only add
slightly more calculations for Lemma 3.3 when $X_{n}$ is defined in (\ref{a0}). (The notations are the same as those given in Remark 3.2 in Lin
and Liu (2009).) Note that
\begin{align}\label{ab1.1}
|r(u)|=|\ep X_{0}X_{u}|=|\sum_{j\in Z}\ep \mathcal{P}_{j}(X_{0})\mathcal{P}_{j}(X_{u})|\leq \sum_{j\in Z}\theta_{j,2}\theta_{u+j,2}
\end{align}
and
\begin{equation}\label{ab1}
\sum_{u\geq n}|r(u)|\leq C\Theta_{[n/2],2}.
\end{equation}
{\bf Lemma 3.3 for $X_{n}$ in (\ref{a0}).}{\em  ~~Suppose that $\ep X_{0}=0$, $\ep X^{2}_{0}<\infty$ and $\Theta_{0,2}<\infty$. Then
\\
(i).
\begin{equation}
\max_{1\leq j\leq q}\Big{|}\frac{\ep S^{2}_{n,j,1}}{\pi nf(\omega_{j})}-1\Big{|}\leq Cn^{-1}\sum_{k=0}^{n}\Theta_{k,2}.
\end{equation}
(ii).
\begin{equation}
\max_{1\leq j\leq q}\Big{|}\frac{\ep S^{2}_{n,j,2}}{\pi nf(\omega_{j})}-1\Big{|}\leq Cn^{-1}\sum_{k=0}^{n}\Theta_{k,2}.
\end{equation}
(iii).   $\max_{1\leq i, j\leq q}|\ep S_{n,i,1}S_{n,j,2}|\leq C\sum_{k=0}^{n}\Theta_{k,2}$ and $\max_{1\leq i\neq j\leq q}|\ep
S_{n,i,l}S_{n,j,l}|\leq C\sum_{k=0}^{n}\Theta_{k,2}$ for $l=1,2$.}

\begin{proof} We only prove {\em (i)}. By (\ref{ab1.1}), (\ref{ab1}) and the proofs in Lin and Liu (2009),
 we have
\begin{align}\label{abcde}
\Big{|}\frac{\ep S^{2}_{n,j,1}}{\pi nf(\omega_{j})}-1\Big{|}&\leq C\sum_{k=n}^{\infty}|r(k)|+Cn^{-1}\sum_{k=1}^{n-1}k|r(k)|\non &\leq
C\Theta_{[n/2],2}+Cn^{-1}\sum_{j=0}^{\infty}\theta_{j,2}\sum_{k=1}^{n}k(\Theta_{k+j,2}-\Theta_{k+j+1,2})\non
&\quad+Cn^{-1}\sum_{j=-n}^{-1}\theta_{j,2}\sum_{k=-j}^{n}k(\Theta_{k+j,2}-\Theta_{k+j+1,2})\non
&\quad+Cn^{-1}\sum_{j=-n}^{-1}\theta_{j,2}\sum_{k=1}^{-j-1}k(\Theta_{|k+j|,2}-\Theta_{|k+j|+1,2})\non
&\quad+Cn^{-1}\sum_{j=-\infty}^{-n-1}\theta_{j,2}\sum_{k=1}^{ n}k\theta_{k+j,2}.
\end{align}
Since $\sum_{k=-j}^{n}k(\Theta_{k+j,2}-\Theta_{k+j+1,2})\leq \sum_{k=0}^{n}\Theta_{k,2}-j\Theta_{0,2}$ for $-n\leq j\leq 0$,
\begin{align*}
&n^{-1}\sum_{j=-n}^{-1}\theta_{j,2}\sum_{k=-j}^{n}k(\Theta_{k+j,2}-\Theta_{k+j+1,2})\\
&\leq
Cn^{-1}\sum_{k=0}^{n}\Theta_{k,2}+Cn^{-1}\sum_{j=-n}^{-1}(-j)\theta_{j,2}\\
&\leq Cn^{-1}\sum_{k=0}^{n}\Theta_{k,2}.
\end{align*}
Similarly, we can show that the other terms in (\ref{abcde}) have the same bound $ Cn^{-1}\sum_{k=0}^{n}\Theta_{k,2}$.
 The proof of the lemma is complete.
\end{proof}

\section{Proof of Lemma 4.2} \setcounter{equation}{0}
In this section, we prove Lemma 4.2 in Lin and Liu (2009).

{\noindent\bf Lemma 4.2 in Lin and Liu (2009)}{\em ~ Let $\xi_{n,1},\cdots, \xi_{n,k_{n}}$ be independent random vectors with mean zero
 and values in $\mathbb{R}^{d}$, and $S_{n}=\sum_{i=1}^{k_{n}}X_{n,i}$. Assume that $|\xi_{n,k}|\leq
 c_{n}B^{1/2}_{n}$, $1\leq k\leq k_{n}$, for some $c_{n}\rightarrow 0$,
 $B_{n}\rightarrow\infty$ and
 \begin{equation*}
 \Big{|}B^{-1}_{n}\Cov(\xi_{n,1}+\cdots+\xi_{n,k_{n}})-I_{d}\Big{|}\leq C_{0}c^{2}_{n},
 \end{equation*}
 where $I_{d}$ is a $d\times d$ identity matrix and $C_{0}$ is a positive constant. Suppose that
$
 \beta_{n}:=B^{-3/2}_{n}\sum_{k=1}^{k_{n}}\ep|\xi_{n,k}|^{3}\rightarrow
 0.
 $
 Then for all $n\geq n_{0}$ ($n_{0}$ is given below)
 \begin{align*}
& |\pr(|S_{n}|_{d}\geq x)-\pr(|N|_{d}\geq
x/B^{1/2}_{n})|\\
&\leq o(1)\pr(|N|_{d}\geq x/B^{1/2}_{n})+C_{d}\Big{(}\exp\Big{(}-\frac{\delta^{2}_{n}\min(c^{-2}_{n}, \beta^{-2/3}_{n})}{8d}\Big{)}
+\exp\Big{(}\frac{C^{-1}_{d}c^{2}_{n}}{\beta^{2}_{n}\log \beta_{n}}\Big{)}\Big{)},
 \end{align*}
uniformly for $x\in[B^{1/2}_{n}, \delta_{n}\min(c^{-1}_{n}, \beta^{-1/3}_{n})B^{1/2}_{n}]$, with any $\delta_{n}\rightarrow 0$ and
$\delta_{n}\min(c^{-1}_{n}, \beta^{-1/3}_{n})\rightarrow \infty$. $N$ is a centered normal random vector with covariance matrix $I_{d}$.
$|\cdot|_{d}$ denotes the
  $d$-dimensional Euclidean norm or $|z|_{d}=\min\{(x^{2}_{i}+y^{2}_{i})^{1/2}: 1\leq i\leq
  d/2\}$, $z=(x_{1}, y_{1},\cdots, x_{d/2},y_{d/2})$ (we assume $d$ is even in this
  case).  $o(1)$ is bounded by $A_{n}:=A(\delta_{n}+\beta_{n}+c_{n})$, $A$ is
  a positive constant  depending only on $d$.
  \begin{eqnarray*}
  &&n_{0}=\min\Big{\{}n: \forall k\geq n,~c^{2}_{k}\leq \frac{\min(C^{-1}_{0},8^{-1})}{2},~\delta_{k}\leq \min((100d)^{-1}c_{17},1,C^{-1}_{d}C^{-2}_{0},C^{-1}_{d}),\cr
&&\qquad\qquad\qquad\qquad\qquad\qquad\qquad\qquad\qquad\qquad\qquad\qquad\qquad\beta_{k}\leq
  \frac{\min(c^{-1}_{22},1)}{32d}\Big{\}},
  \end{eqnarray*}
  where $c_{17}$ and $c_{22}$ (given in  Corollary 1 in Einmahl \cite{einmahl}) are positive constants depending only on $d$, and $C_{d}$ (given in the
  proof) is a positive constant depending only on $d$.
}

\begin{proof}  Set
$\Sigma_{n}=\Cov(\xi_{n,1}+\cdots+\xi_{n,k_{n}})$, $\xi^{'}_{n,k}=B_{n}^{1/2}\Sigma^{-1/2}_{n}\xi_{n,k}$, $1\leq k\leq k_{n}$, and
$S^{'}_{n}=\sum_{k=1}^{k_{n}}\xi^{'}_{n,k}$. Then
\begin{equation*}
\Cov(\xi^{'}_{n,1}+\cdots+\xi^{'}_{n,k_{n}})=B_{n}I_{d}.
\end{equation*}
Note that for $n\geq n_{0}$, $|\xi^{'}_{n,k}|\leq 2c_{n}B^{1/2}_{n}$ for $1\leq k\leq k_{n}$. We now use Corollary 1(b) in Einmahl \cite{einmahl}
to prove the lemma. Taking $\alpha=(100dc_{n}B^{1/2}_{n})^{-1}$ in that corollary, it can be checked that for $n\geq n_{0}$,
\begin{equation*}
\alpha\sum_{k=1}^{k_{n}}\ep|\xi^{'}_{n,k}|^{3}\exp(\alpha|\xi^{'}_{n,k}|)\leq B_{n}.
\end{equation*}
Write $ \beta^{'}_{n}=B^{-3/2}_{n}\sum_{k=1}^{k_{n}}\ep|\xi^{'}_{n,k}|^{3}\exp(\alpha|\xi^{'}_{n,k}|) $ and it holds that $\beta^{'}_{n}\leq
16d\beta_{n}$ for $n\geq n_{0}$. Let $\eta_{1},\cdots,\eta_{k_{n}}$ be independent $N(0, \sigma^{2}\Cov(\xi^{'}_{n,k}))$ random vectors, which are
also independent of $\{\xi^{'}_{n,k}\}$ and $0<\sigma^{2}\leq 1$.  Take $\sigma^{2}= -c_{22}(16d)^{2}\beta^{2}_{n}\log (16d\beta_{n})$ ($c_{22}$
is defined in (\ref{ein})) and set $M_{n}=B^{-1/2}_{n}\Sigma^{1/2}_{n}-I_{d}$. Then $|M_{n}|^{2}=|M^{2}_{n}|\leq
|B^{-1}_{n}\Sigma_{n}-I_{d}|^{2}\leq C^{2}_{0}c^{4}_{n}$. Since $|x+y|_{d}\leq |x|_{d}+|y|$ for $x, y\in R^{d}$, we have
\begin{align}\label{a7}
&\pr(|S_{n}|_{d}\geq x)= \pr(|M_{n}S^{'}_{n}+S^{'}_{n}|_{d}\geq x)\non &\leq \pr\Big{(}\Big{|}S^{'}_{n}+\sum_{k=1}^{k_{n}}\eta_{k}\Big{|}_{d}\geq
x-c_{n}B^{1/2}_{n}\Big{)}+\pr\Big{(}\Big{|}\sum_{k=1}^{k_{n}}\eta_{k}\Big{|}\geq c_{n}B^{1/2}_{n}/2\Big{)}\non &\quad+\pr\Big{(}
|M_{n}S^{'}_{n}|\geq c_{n}B^{1/2}_{n}/2\Big{)}\non &\leq \pr\Big{(}\Big{|}S^{'}_{n}+\sum_{k=1}^{k_{n}}\eta_{k}\Big{|}_{d}\geq
x-c_{n}B^{1/2}_{n}\Big{)}+C_{d}\exp\Big{(}\frac{C^{-1}_{d}c^{2}_{n}}{\beta^{2}_{n}\log \beta_{n}}\Big{)}+
C_{d}\exp\Big{(}-C^{-1}_{d}\min(C^{-2}_{0},1)c^{-2}_{n}\Big{)},
\end{align}
where the last inequality follows from the exponential inequality (cf. Lemma 1.6 in Ledoux and Talagrand \cite{ledoux}) and  $C_{d}$ is a positive
constant depending only on $d$.

  Corollary
1(b) (in combination with the Remark on page 32) in Einmahl \cite{einmahl} implies that,
 if
\begin{equation}\label{ein}
|x|\leq c_{17}\alpha B^{1/2}_{n},\quad 1\geq\sigma^{2}\geq -c_{22}\beta^{'2}_{n}\log \beta^{'}_{n}\quad\mbox{and }B_{n}\geq c_{18}\alpha^{-2},
\end{equation}
where $c_{17}$, $c_{22}$ and $c_{18}$ (given in Einmahl \cite{einmahl}) are constants depending only on $d$, then
\begin{equation*}
p_{n}(x)=\varphi_{(1+\sigma^{2})I_{d}}(x)\exp(T_{n}(x)), \quad \mbox{with } |T_{n}(x)|\leq c_{19}\beta^{'}_{n}(|x|^{3}+1),
\end{equation*}
where $p_{n}(x)$ is the density of $B^{-1/2}_{n}\sum_{k=1}^{k_{n}}(\xi^{'}_{n,k}+\eta_{k})$, $\varphi_{M}$ is the density of a $d$-dimensional
centered Gaussian vector with covariance matrix $M$, and $c_{19}$ is a constant only depending on $d$. Letting
$t_{n}=\delta_{n}\min\{\beta^{-1/3}_{n}, c^{-1}_{n}\}$ and noting that $|x|\leq t_{n}$ implies $|x|\leq c_{17}\alpha B^{1/2}_{n}$ for $n\geq
n_{0}$, we have
\begin{align}\label{a8}
&\pr\Big{(}\Big{|}S^{'}_{n}+\sum_{k=1}^{k_{n}}\eta_{k}\Big{|}_{d}\geq x-c_{n}B^{1/2}_{n}\Big{)}\non &=\int_{|y|_{d}\geq
x/B^{1/2}_{n}-c_{n}}p_{n}(y)dy\non &=\int_{|y|_{d}\geq x/B^{1/2}_{n}-c_{n},|y|\leq t_{n}}\varphi_{(1+\sigma^{2})I_{d}}(y)\exp(T_{n}(y))dy\non
&\quad+\int_{|y|_{d}\geq x/B^{1/2}_{n}-c_{n},|y|> t_{n}}p_{n}(y)dy.
\end{align}
The first term on the right hand side of (\ref{a8}) is
\begin{align}\label{a9}
&\int_{|y|_{d}\geq x/B^{1/2}_{n}-c_{n},|y|\leq t_{n}}\varphi_{(1+\sigma^{2})I_{d}}(y)\exp(T_{n}(y))dy\non
 &\leq
(1+A_{n})\int_{|y|_{d}\geq x/B^{1/2}_{n}-c_{n},|y|\leq t_{n}}\varphi_{(1+\sigma^{2})I_{d}}(y)dy\non &=(1+A_{n})\int_{|y|_{d}\geq
x/B^{1/2}_{n}-c_{n}}\varphi_{(1+\sigma^{2})I_{d}}(y)dy\non &\quad-(1+A_{n})\int_{|y|_{d}\geq x/B^{1/2}_{n}-c_{n},|y|\geq
t_{n}}\varphi_{(1+\sigma^{2})I_{d}}(y)dy.
\end{align}
For $x\in[B^{1/2}_{n}, \delta_{n}\min(c^{-1}_{n}, \beta^{-1/3}_{n})B^{1/2}_{n}]$, we have $c_{n}x/B^{1/2}_{n}\leq \delta_{n}$. This together with
some elementary calculations implies that
\begin{align}\label{a10}
&\int_{|y|_{d}\geq x/B^{1/2}_{n}-c_{n}}\varphi_{(1+\sigma^{2})I_{d}}(y)dy\non &\leq \pr(|N|_{d}\geq x/B^{1/2}_{n}-2c_{n})+\pr(\sigma|N|\geq
c_{n})\non &\leq(1+A_{n})\pr(|N|_{d}\geq x/B^{1/2}_{n})+C_{d}\exp\Big{(}\frac{C^{-1}_{d}c^{2}_{n}}{\beta^{2}_{n}\log \beta_{n}}\Big{)}
\end{align}
and
\begin{align}\label{a11}
\int_{|y|_{d}\geq x/B^{1/2}_{n}-c_{n},|y|\geq t_{n}}\varphi_{(1+\sigma^{2})I_{d}}(y)dy\leq C_{d}\exp\Big{(}-\frac{\delta^{2}_{n}\min(c^{-2}_{n},
\beta^{-2/3}_{n})}{4}\Big{)}.
\end{align}
For the second term in (\ref{a8}), we shall use again Lemma 1.6 in Ledoux and Talagrand \cite{ledoux}, and it follows that
\begin{align}\label{a12}
\int_{|y|_{d}\geq x/B^{1/2}_{n}-c_{n},|y|> t_{n}}p_{n}(y)dy &\leq\pr\Big{(}|\sum_{k=1}^{k_{n}}\xi^{'}_{n,k}|\geq
9t_{n}B^{1/2}_{n}/10\Big{)}+\pr(\sigma|N|\geq t_{n}/10)\non &\leq C_{d}\exp\Big{(}-\frac{\delta^{2}_{n}\min(c^{-2}_{n},
\beta^{-2/3}_{n})}{8d}\Big{)}.
\end{align}
Finally, combining (\ref{a7})-(\ref{a12}) gives
\begin{align*}
 \pr(|S_{n}|_{d}\geq x) &\leq (1+A_{n})(\pr(|N|_{d}\geq
x/B^{1/2}_{n}))\\
&+C_{d}\Big{(}\exp\Big{(}-\frac{\delta^{2}_{n}\min(c^{-2}_{n}, \beta^{-2/3}_{n})}{8d}\Big{)} +\exp\Big{(}\frac{C^{-1}_{d}c^{2}_{n}}{\beta^{2}_{n}\log
\beta_{n}}\Big{)}\Big{)}.
 \end{align*}
Similarly, we can show that
\begin{align*}
\pr(|S_{n}|_{d}\geq x) &\geq (1-A_{n})(\pr(|N|_{d}\geq
x/B^{1/2}_{n}))\\
&-C_{d}\Big{(}\exp\Big{(}-\frac{\delta^{2}_{n}\min(c^{-2}_{n}, \beta^{-2/3}_{n})}{8d}\Big{)} +\exp\Big{(}\frac{C^{-1}_{d}c^{2}_{n}}{\beta^{2}_{n}\log
\beta_{n}}\Big{)}\Big{)}.
 \end{align*}
The desired result now follows.
\end{proof}


\begin{thebibliography}{99}

\bibitem{an} An, H. Z., Chen, Z. G. and Hannan, E. J. (1983). The maximum of the periodogram.
{\em Journal of Multivariate Analysis, }{\bf 13: } 383-400.

\bibitem{brock}  Brockwell, P. J. and Davis, R. A. (1998). Time
Series: Theory and Methods, 2nd ed. New York: Springer-Verlag.


\bibitem{chiu}  Chiu, S. T. (1989). Detecting periodic components in a white Gaussian time series. {\em
J. R. Statist. Soc.  } B. {\bf 51:} 249-259.

\bibitem{davis} Davis, R. A. and Mikosch, T.
 (1999). The maximum of the periodogram of a non-Gaussian sequence. {\em
 Ann. Probab.} {\bf 27: } 522-536.

\bibitem{einmahl} Einmahl, U. (1989). Extensions of results of Koml\'{o}s, Major, and Tusn\'{a}dy to the
multivariate case. {\em J. Mult. Anal.,} {\bf 28:} 20-68.


\bibitem{fisher} Fisher, R. A. (1929). Tests of significance in harmonic
analysis.   {\em Proc. Roy. Statist. Soc., Ser. }  {\bf 125:} 54-59.

\bibitem{freed} Freedman, D. A. (1975). On tail probabilities for
martingales. {\em Ann. Probab.}, {\bf 3: } 100-118.

\bibitem{lewis} Lewis, T. and Fieller, N. R. J. (1979). A recursive
algorithm for null distributions for outliers: I. Gamma samples. {\em Technometrics.} {\bf 21: } 371-376.






\bibitem{hannan1}  Hannan, E. J. (1961). Testing for a jump in the spectral function. {\em
J. R. Statist. Soc.  } B. {\bf 23:} 394-404.

\bibitem{hannan2}  Hannan, E. J. (1970). Multiple Time Series, pp.
463-475. New York: Wiley.

\bibitem{hsing} Hsing, T and Wu, W. B. (2004). On weighted U-statistics for stationary
processes. {\em Ann. Probab., } 32, 1600-1631.

\bibitem{lin} Lin, Z. Y. and Liu, W. D. (2008). Supplementary material for "On maxima of periodograms  of stationary
processes". http://www.math.zju.edu.cn/stat/linliu.

pdf; http://xxx.arxiv.org/abs/0801.1357.

\bibitem{nagaev} Nagaev, S. V. (1979). Large deviations of independent random variables. {\em Ann.
Probab.} 7, 745-789.

\bibitem{priestley} Priestley, M. B. (1981). Spectral Analy sis and Time Series. Academic Press, London.

\bibitem{shaox} Shao, X. and Wu, W. B. (2007). Asymptotic spectral theory for nonlinear time
series. {\em Ann. Statist.} {\bf 35 (4):} 1773-1801.

\bibitem{shim} Shimshoni, M. (1971). On Fisher's test of
significance in harmonic analysis. {\em Geophys. J. R. Astronom. Soc. } {\bf 23: } 373-377.

\bibitem{walker} Walker, A. M. (1965). Some asymptotic results for
the periodogram of a stationary time series. {\em J. Aust. Math. Soc.} {\bf 5:} 107-128.

\bibitem{wood} Woodroofe, M. and Van Ness, J. W. (1967). The maximum
deviation of sample spectral densities. {\em Ann. Math. Statist.} {\bf 38: } 1558-1569.


\bibitem{wu} Wu, W. B. (2005) Nonlinear system theory: another look at
dependence. {\em Proc. Natl Acad. Sci. USA. } {\bf 102(40):} 14150-14154.

\bibitem{wu1} Wu, W. B. (2007). Strong invariance principles for
dependent random variables. {\em Ann. Probab.,} {\bf 35: } 2294-2320.

\bibitem{wu3} Wu, W. B. and Shao, X. (2004). Limit theorems for iterated random
functions. {\em Journal of Applied Probability, } {\bf 41:} 425-436.





















\end{thebibliography}

\begin{thebibliography}{99}
\bibitem{einmahl} Einmahl, U. (1989). Extensions of results of Koml\'{o}s, Major, and Tusn\'{a}dy to the
multivariate case. {\em J. Mult. Anal.,} {\bf 28:} 20-68.

\bibitem{lin} Lin, Z.Y. and Liu, W.D. (2009). On maxima of periodograms
of stationary processes. {\em Ann. Statist.} {\bf 37: } 2676-2695.

\bibitem{ledoux} Ledoux, M. and Talagrand, M. (1991). Probability
in Banach Spaces. Springer, Berlin.
\end{thebibliography}
\end{document}